\documentclass{amsart}

\usepackage{amsfonts}
\usepackage{CJK}
\usepackage{fancyhdr}
\usepackage{graphicx}
\usepackage[dvips,usenames]{color}
\usepackage{titletoc}
\usepackage{latexsym}
\usepackage{amssymb}
\usepackage{multicol}
\usepackage{graphics}
\usepackage{subfigure}
\usepackage{indentfirst}
\usepackage{cases}
\usepackage{curves}

\newtheorem{theorem}{Theorem}[section]
\newtheorem{corollary}[theorem]{Corollary}
\newtheorem{lemma}[theorem]{Lemma}
\newtheorem{definition}[theorem]{Definition}

\numberwithin{equation}{section}

\def\x#1{(\ref{#1})}
\def\R{{\Bbb R}}

\def\Z{{\Bbb Z}}

\newcommand{\cF}{\mathcal{F}}

\newcommand{\cS}{\mathcal{S}}
\newcommand{\cU}{\mathcal{U}}
\newcommand{\cW}{\mathcal{W}}
\newcommand{\cI}{\mathcal{I}}

\newcommand{\ned}{~nonuniform exponential dichotomy}

\newcommand{\nds}{nonuniform dichotomy spectrum~}

\def\bb{\begin}
\def\bc{\begin{center}}
\def\ec{\end{center}}
\def\ba{\begin{array}}
\def\ea{\end{array}}
\def\be{\begin{equation}}
\def\ee{\end{equation}}
\def\bea{\begin{eqnarray}}
\def\eea{\end{eqnarray}}
\def\beaa{\begin{eqnarray*}}
\def\eeaa{\end{eqnarray*}}

\def\ben{\begin{enumerate}}
\def\een{\end{enumerate}}
\def\hh{\!\!\!\!}

\def\EQ{\hh & = & \hh}

\def\LE{\hh & \le & \hh}

\def\nn{\nonumber}
\def\oo{\infty}
\def\ifl{\iffalse}
\def\lb{\label}
\def\prf{\mbox{\bf Proof.~}}

\DeclareMathOperator{\im}{im}

\title[]{Nonuniform dichotomy spectrum and reducibility for nonautonomous equations}

\author[ ]
{Jifeng Chu$^{1,2}$, \quad Fang-Fang Liao $^1$,\quad Stefan
Siegmund$^2$,\quad Yonghui Xia$^3$, Weinian Zhang $^4$}
\address{$^1$ Department of Mathematics, College of Science,
Hohai University, Nanjing 210098, China}
\address{$^2$ Center for Dynamics \& Institute for Analysis, Department of Mathematics, TU Dresden, Germany}
\address{$^3$ College of Mathematics, Physics and Information Engineering, Zhejiang
Normal University, Jinhua, China}
\address{$^4$ College of Mathematics, Sichuan University, Chengdu, China}
\email{jifengchu@126.com (J. Chu)}
\email{liaofangfang8178@sina.com (F. Liao)}
\email{stefan.siegmund@tu-dresden.de (S. Siegmund)}
\email{yhxia@zjnu.cn; xiadoc@163.com (Y. Xia)}
\email{wnzhang@scu.edu.cn  (W. Zhang)}
\thanks{Jifeng Chu was supported by the National Natural Science Foundation
of China (Grant No. 11171090, No. 11271078), China Postdoctoral Science
Foundation funded project (Grant No.2012T50431) and the Alexander von Humboldt
Foundation of Germany. Yonghui Xia was supported by the Natural Science Foundation
of China (Grant No.11271333).}


\begin{document}

\begin{abstract} For nonautonomous linear differential equations with nonuniform
hyperbolicity, we introduce a definition for nonuniform dichotomy spectrum, which can
be seen as a generalization of Sacker-Sell spectrum. We prove a spectral theorem and
use the spectral theorem to prove a reducibility result.
\end{abstract}

\maketitle


\section{\bf Introduction}
\setcounter{equation}{0}

Let $ {\mathfrak{L}}_{loc}^1= {\mathfrak{L}}_{loc}^1(\mathbb{R},\mathbb{R}^{N\times N}),N\in \mathbb{N}$,
be the space of locally integrable matrix functions. Given $A \in {\mathfrak{L}}_{loc}^1$,
we consider the following nonautonomous linear differential equation
\begin{equation}x'=A(t)x.\label{1.1}\end{equation}
Let $\Phi:\R\times\R\rightarrow\R^{N\times N},~(t,s)\mapsto \Phi(t,s)$ denote the associated evolution operator of (\ref{1.1}),
i.e., $\Phi(t,s)x(s)=x(t)$ for every $t,s \in \mathbb{R}$, where $x$ is any solution of (\ref{1.1}).
Clearly, $\Phi(t,\tau)\Phi(\tau,s)=\Phi(t,s)$, $t,\tau,s\in \mathbb{R}$.

The classical notion of exponential dichotomy introduced by Perron in \cite{per}
plays an important role in the study of dynamical behaviors of \x{1.1}, particularly
in what concerns the study of stable and unstable invariant manifolds, and therefore
has attracted much attention (see, for example, \cite{cl94, cl95, lat96, Pal84, pot09, Rod88, ss78, ss94, Zh95})
during the last few decades.
We also refer to the books \cite{cl, cop78, ms} for details and
further references related to exponential dichotomies.
On the other hand, as Barreira and Valls mentioned in \cite{bv08},
the classical notion of exponential dichotomy substantially restricts some
dynamics and it is important to look for more general types of hyperbolic
behaviors. During the last several years, inspired by both the classical notion of exponential dichotomy and
the notion of nonuniformly hyperbolic trajectory introduced
by Pesin in \cite{bp02, bp07}, Barreira and Valls introduced the concept
of nonuniform exponential dichotomy and investigated some related problems
\cite{bv05, bv06, bv07, bv10}. In particular, they discussed the existence and the
smoothness of invariant manifolds for nonautonomous differential equations,
a version of the Grobman-Hartman theorem, the existence
of center manifolds and the theory of Lyapunov regularity. A more general nonuniform
exponential dichotomy has been considered in \cite{bcv-1, bcv-2, chu}, which admits
different growth rates in the uniform and nonuniform parts.
Barreira and Valls explained in \cite{bp02, bv08} that, from the point of view of ergodic theory, almost all linear
variational equations have a nonuniform exponential behavior.

Based on the study of classical exponential dichotomy, the dichotomy spectral theory
was introduced by Sacker and Sell in \cite{ss78}.
The dichotomy spectrum is an important object in the theory of dynamical systems
because the spectral intervals, together with the spectral manifolds,
completely describe the dynamical skeleton of a linear system. A spectral theory based on
finite-time hyperbolicity has been studied in \cite{bds, dkns, dps}.
Some other related results can be seen from
\cite{am, amz, as01, as02, cl94, ns, pot09, pot12, s-jdde}.
The dichotomy spectral theory was applied in \cite{s-jlms, s-jde} to give
block diagonalization and normal forms for nonautonomous differential equations.


In this paper we investigate the dichotomy spectrum in the setting of nonuniform exponential dichotomies,
called {\em nonuniform dichotomy spectrum}. We show the topological structure of the nonuniform dichotomy spectrum
and give the corresponding decomposition in spectral manifolds.
At last, we use the above results on spectrum to prove the reducibility for \x{1.1}, i.e., a kinematical
similarity to a diagonal system in proper blocks.


\section{Nonuniform dichotomy spectrum}
\setcounter{equation}{0}

\noindent Let $\Phi(t,s)$ be the evolution operator of (\ref{1.1}).
An {\em invariant projector} of (\ref{1.1}) is defined to be a function
$P:\mathbb{R}\rightarrow\mathbb{R}^{N\times N}$
of projections $P(t),\,t\in \mathbb{R}$, such that
\[P(t)\Phi(t,s)=\Phi(t,s)P(s),\,\,\,\mbox{for} \,\,\,\,t,s \in \mathbb{R}.\]
Clearly, $P$ is continuous due to the identity $P\equiv \Phi(\cdot,s)P(s)\Phi(s,\cdot)$.

\begin{definition} \label{Def1.1}
We say that equation {\rm(\ref{1.1})} admits a nonuniform exponential
dichotomy on $\mathbb{R}$, if there exist constants $\alpha>0$, $K>0$, $\varepsilon\geq0$ with $\varepsilon<\alpha$
and an invariant projector $P$ such that
\be\lb{np}\|\Phi(t,s)P(s)\|\leq Ke^{-\alpha(t-s)}e^{\varepsilon |s|}, \quad t\geq s,\ee
and
\be\lb{nq} \|\Phi(t,s)Q(s)\|\leq Ke^{\alpha(t-s)}e^{\varepsilon |s|},\quad t\leq s,\ee
where $Q(t)={\rm Id}-P(t)$ is the complementary projection.
\end{definition}
When one can take $\varepsilon=0$ in \x{np}-\x{nq}, we say that equation (\ref{1.1})
admits a (uniform) exponential dichotomy, and thus a classical exponential dichotomy is a particular
case of a nonuniform one. As illustrated in \cite{bv08}, in most cases, the nonuniform part $e^{\varepsilon |s|}$ in \x{np}-\x{nq} can not be removed.
In particular, Barreira and Valls have proved that in finite-dimensional spaces essentially any linear equation with nonzero Lyapunov
exponents admits a nonuniform exponential dichotomy, and as a consequence of Oseledets multiplicative ergodic theorem \cite{o},
the nonuniformity of most equations is very small. See \cite{bv08} for the details. We remark that \cite[Theorem 1.4.2]{bp02} indicates that the condition
$\varepsilon<\alpha$ is reasonable,
which means that the nonuniform parts are small.

For example, if $\lambda>3a>0$, then the linear equation in $\mathbb{R}^2$ given by \[u'=(-\lambda-at\sin t)u,\,\,\,v'=(\lambda+at\sin t)v,\]
admits a nonuniform exponential dichotomy, but it does not admit a (uniform) exponential dichotomy.

\begin{lemma}\label{lemma2.1} The projector of equation {\rm\x{1.1}}
can be chose as \[\tilde{P}= \left(\begin{array}{lll}
I_{N_1} & 0\\
0 &  0_{N_2}\end{array}\right).\] with
$N_1 = {\rm dim\, im} \tilde{P}$ and $N_2 = {\rm dim\, ker} \tilde{P}$, and
the fundamental matrix $X(t)$ can be chosen appropriately such that the
estimates {\rm\x{np}-\x{nq}} can be rewritten as
\be\lb{ntp}\|X(t)\tilde{P}X^{-1}(s)\|\leq Ke^{-\alpha(t-s)}e^{\varepsilon |s|},\quad t\geq s,\ee
and
\be\lb{ntq}\|X(t)\tilde{Q}X^{-1}(s)\|\leq Ke^{\alpha(t-s)}e^{\varepsilon |s|}, \quad t\leq s,\ee
where $\tilde{Q}={\rm Id}-\tilde{P}.$
\end{lemma}

\noindent{\bf Proof.}\,
Let $\tau\in \mathbb{R}$ be arbitrarily chosen but fixed. Then there exists a non-singular matrix $T\in \mathbb{R}^{N\times N}$ such that
\[TP(\tau)T^{-1}= \left(\begin{array}{lll}
I_{N_1} & 0\\
0 &  0_{N_2}\end{array}\right).\]
For the evolution operator $\Phi(t,\tau)$ of (\ref{1.1}), we define
\[X(t) := \Phi(t,\tau)T^{-1},\quad t\in\R,\]
and
\[\tilde{P}:=  \left(\begin{array}{lll}
I_{N_2} & 0\\0 &  0_{N_2}\end{array}\right)= T P(\tau) T^{-1}.\]
Then
\bea
\|X(t)\tilde{P}X^{-1}(s)\|\EQ\|\Phi(t,\tau)T^{-1} \tilde{P}T \Phi^{-1}(s,\tau)\|\nn\\
\EQ\|\Phi(t,\tau)P(\tau)\Phi^{-1}(s,\tau)\|.
\label{2.3}
\eea
On the other hand, one has
\begin{equation}\begin{array}{lll}
\|\Phi(t,s)P(s)\|&=&\|\Phi(t,\tau)\Phi(\tau,s)P(s)\|\\
&=&\|\Phi(t,\tau)P(\tau)\Phi(\tau,s)\|\\
&=&\|\Phi(t,\tau)P(\tau)\Phi^{-1}(s,\tau)\|.
\end{array}\label{2.4}\end{equation}
It follows from (\ref{2.3})-(\ref{2.4}) that \x{np}-\x{nq} can be
rewritten in the equivalent form \x{ntp}-\x{ntq}. $\Box$

For fixed $\gamma \in \mathbb{R}$,
consider the shifted system
\be\dot{x} = [A(t) - \gamma I]x, \tag*{$(2.8)_\gamma$}\ee
\refstepcounter{equation}
which has the evolution operator
\[
\Phi_\gamma(t,s) := e^{-\gamma (t-s)}\Phi(t,s).
\]
If (2.8)$_\gamma$ admits a nonuniform exponential dichotomy, then its
invariant projector $P(t)$ is also invariant for \x{1.1}. The dichotomy estimates are
equivalent to
\be\lb{Px}\|\Phi(t,s)P(s)\|\leq K e^{(\gamma-\alpha)(t-s)}e^{\varepsilon |s|}, \quad t \geq s\ee
and
\be\lb{Qx}\|\Phi(t,s)Q(s)\|\leq K e^{(\gamma + \alpha)(t-s)}e^{\varepsilon |s|}, \quad t \leq s.\ee
By Lemma \ref{lemma2.1}, equivalently,
\[X_\gamma(t) := e^{-\gamma t}X(t)=e^{-\gamma t}\Phi(t,\tau)T^{-1}\]
is the fundamental matrix of the shifted system (2.8)$_\gamma$, and its invariant projection is \[\tilde{P}= \left(\begin{array}{lll}
I_{N_1} & 0\\0 &  0\end{array}\right).\]
The corresponding estimates are equivalent to
\be\lb{xnp}\|X_\gamma(t)\tilde{P}X_\gamma^{-1}(s)\|\leq Ke^{-\alpha(t-s)}e^{\varepsilon |s|},\quad t\geq s,\ee
and
\be\lb{xnq}\|X_\gamma(t)\tilde{Q}X_\gamma^{-1}(s)\|\leq Ke^{\alpha(t-s)}e^{\varepsilon |s|},\quad t\leq s.\ee
We will use the estimates \x{Px}-\x{Qx} as well as the equivalent formulation \x{xnp}-\x{xnq}.

\begin{definition} \label{Def3.1}
The nonuniform dichotomy spectrum of {\rm(\ref{1.1})} is the set
\[\Sigma_{NED}(A) = \{ \gamma\in\R:(2.8)_\gamma \mbox{ admits no nonuniform exponential dichotomy} \},\]
and the \emph{resolvent set} $\rho_{NED}(A) =  \mathbb{R} \setminus \Sigma_{NED}(A)$ is
its complement.
\end{definition}

Let $\Sigma_{ED}(A)$ denote the classical dichotomy spectrum of \x{1.1}. Obviously,
$\Sigma_{NED}(A)\subset \Sigma_{ED}(A)$.
For $\gamma \in\rho_{NED}(A)$, define
\[\cS_\gamma := \left\{(\tau,\xi) \in \R \times \R^N \,:\,
\sup_{t\geq 0}\{\|\Phi(t,\tau)\xi\|e^{-\gamma t} \} e^{- \varepsilon\tau}<\oo\right\},\]
and
\[\cU_\gamma := \left\{(\tau,\xi) \in \R \times \R^N \,:\,
\sup_{t\leq 0}\{\|\Phi(t,\tau)\xi\|e^{-\gamma t} \} e^{- \varepsilon\tau}<\oo\right\}.\]
One may readily verify that $\cS_\gamma$ and $\cU_\gamma$ are linear integral manifold of \x{1.1}.
As defined in \cite{s-jdde}),
a nonempty set $\mathcal{W} \subset \R \times \R^N$
is a {\it linear integral manifold of} {\rm(\ref{1.1})} if
(a) it is {\it invariant}, i.e., $(\tau,\xi) \in \mathcal{W} \;\Rightarrow\; (t,\Phi(t,\tau)\xi) \in \mathcal{W}$
for all $t \in \R$, (b) for every $\tau \in \R,$ the {\it fiber}
$\mathcal{W}(\tau) = \{ \xi \in \R^N \,:\, (\tau,\xi) \in \mathcal{W} \}$
is a linear subspace of $\R^N$.

At first glance, $\cS_\gamma$ and $\cU_\gamma$ are not well defined
because they seem to depend on the constant $\varepsilon$, which may
not be unique in \x{np}-\x{nq}. However, the following
result ensures that $\cS_\gamma$ and $\cU_\gamma$ are well defined in the setting of a \ned~and they
do not depend on the choice of the constant $\varepsilon$. First we recall that
the invariant projector $P$ is unique for \x{1.1} following the arguments in \cite[Chapter 2]{cop78}.
Although the arguments in \cite{cop78} are done in the setting of exponential dichotomies,
it is not difficult to verify that they are also applicable to the case of nonuniform exponential dichotomies.

\begin{lemma}  \label{lem2.1}
Assume that {\rm (2.8)$_\gamma$} admits a \ned ~with invariant projector $P$ for
$\gamma \in\rho_{NED}(A).$ Then
\[\cS_\gamma = \im P\;,\quad\cU_\gamma = \ker P\quad \text{and} \quad
\cS_\gamma \oplus \cU_\gamma = \R \times \R^N\;. \]
\end{lemma}

\noindent\prf{We show only $\cS_\gamma = \im P$. The fact $\cU_\gamma = \ker P$ is analog and the fact
$\cS_\gamma \oplus \cU_\gamma = \R\times \R^N$ is clear.

First we show $\cS_\gamma\subset\im P.$ Let $\tau\in\R$ and $\xi\in\cS_\gamma(\tau).$ Then there exists
a positive constant $C$ such that
\[\|\Phi(t,\tau)\xi\|\leq C e^{\gamma t}e^{\varepsilon\tau},\quad t\geq \tau.\]
We write $\xi=\xi_1+\xi_2$ with $\xi_1\in {\rm im} P(\tau)$ and
$\xi_2\in {\rm ker} P(\tau).$ We show that $\xi_2=0.$ The invariance of
$P$ implies for $t\in\R$ that we have the equivalence
\[\xi_2=\Phi_\gamma(\tau,t)\Phi_\gamma(t,\tau)Q(\tau)\xi=\Phi_\gamma(\tau,t)Q(t)\Phi_\gamma(t,\tau)\xi.\]
Since (2.8)$_\gamma$ admits a \ned, the following inequality holds
\[\|\Phi_\gamma(\tau,t)Q(t)\|\leq K e^{\alpha(\tau-t)}e^{\varepsilon|t|},\quad t\geq \tau.\]
Thus\beaa\|\xi_2\|\LE Ke^{\alpha(\tau-t)}e^{\varepsilon|t|}\|\Phi_\gamma(t,\tau)\xi\|\\
\LE KCe^{\alpha(\tau-t)}e^{\varepsilon|t|}e^{\gamma t}e^{\varepsilon\tau}e^{-\gamma(t-\tau)}\\
\EQ KCe^{\alpha(\tau-t)}e^{\varepsilon|t|}e^{\varepsilon\tau}e^{\gamma\tau}\\
\LE KCe^{(\alpha-\varepsilon)(\tau-t)}e^{\varepsilon|\tau|}e^{\varepsilon\tau}e^{\gamma\tau}\eeaa
which implies that $\xi_2=0$ by letting $t\rightarrow \oo$, since $\varepsilon<\alpha$.

Next we show ${\rm \im} P\subset \cS_\gamma.$ Let $\tau\in\R$ and
$\xi\in {\rm im} P(\tau),$ i.e., $P(\tau)\xi=\xi.$ The \ned ~implies that
\[\|\Phi_\gamma(t,\tau)\xi\|\leq Ke^{-\alpha(t-\tau)}e^{\varepsilon|\tau|}\|\xi\|\leq Ke^{\varepsilon|\tau|}\|\xi\|,\quad t\geq \tau,\]
since $\alpha>0$, which implies that
\[\|\Phi(t,\tau)\xi\|\leq Ke^{-\gamma(t-\tau)}e^{\varepsilon|\tau|}\|\xi\|,\]
and hence $\xi\in\cS_\gamma(\tau).$} $\Box$

\begin{lemma}\label{lemma3.1}  The resolvent set is open, i.e., for every $\gamma \in \rho_{NED}(A),$ there
exists a constant $\beta =\beta(\gamma)>0$ such that $(\gamma-\beta,
\gamma+\beta) \subset \rho_{NED}(A)$. Furthermore,
\[\cS_\zeta = \cS_\gamma,\quad
\cU_\zeta =  \cU_\gamma\quad \text{for} \quad\zeta \in(\gamma-\beta,\gamma+\beta).\]
\end{lemma}

\noindent{\bf Proof.} Let $\gamma \in \rho_{NED}(A)$. Then (2.8)$_\gamma$
admits a \ned, i.e., the estimates \x{xnp}-\x{xnq} hold
with an invariant projector $\tilde{P}$ and constants $K \geq 0$,
$\alpha > 0$ and $\varepsilon\geq0$. For $\beta:= \alpha / 2>0$ and $\zeta \in (\gamma-\beta,
\gamma+\beta)$ we have
\[X_\zeta(t) = e^{(\gamma - \zeta)t} X_\gamma(t).\] Now $\tilde{P}$ is also an
invariant projector for
\[\dot{x} = [A(t) - \zeta I]x\] and we have the estimates
\[\|X_\zeta(t)\tilde{P}X_\zeta^{-1}(s)\|\leq K e^{(\gamma - \zeta - \alpha)(t-s)}e^{\varepsilon |s|}
\leq  K e^{-\beta (t-s)} e^{\varepsilon |s|},\quad t \geq s,\]
and
\[\|X_\zeta(t)\tilde{P}X_\zeta^{-1}(s)\| \leq K e^{(\gamma - \zeta + \alpha)(t-s)} e^{\varepsilon |s|}
\leq  K e^{\beta (t-s)}e^{\varepsilon |s|}, \quad t \leq s.\]
Hence $\zeta \in \rho_{NED}(A)$ and therefore $\rho_{NED}(A)$ is an open set. $\Box$

\begin{corollary}\label{cor3.1}
$\Sigma_{NED}(A)$ is a closed set.
\end{corollary}

Using the facts proved above, we can obtain the following result, whose proof
is similar to \cite[Lemma 3.2]{s-jdde}, and therefore we omit the proof here.

\begin{lemma}  \label{lemma2.3}
Let $\gamma_1,\gamma_2 \in \rho_{NED}(A)$ with $\gamma_1 <
\gamma_2$. Then $\cF
= \cU_{\gamma_1} \cap \cS_{\gamma_2}$ is a linear integral manifold
which satisfies exactly one of the following two alternatives and
the statements given in each alternative are equivalent:

\begin{tabular}{ll}
\hspace{1.5cm}Alternative I & \hspace{1.5cm}Alternative II
\\[1ex]
{\rm(A)} $\cF = \Z \times \{0\}$. & {\rm(A')} $\cF \not= \Z \times
\{0\}$.
\\[0.5ex]
{\rm(B)} $[\gamma_1,\gamma_2] \subset \rho_{NED}(A)$. & {\rm(B')}
There is a $\zeta \in (\gamma_1,\gamma_2) \cap \Sigma_{NED}(A)$.
\\[0.5ex]
{\rm(C)} $\cS_{\gamma_1} = \cS_{\gamma_2}$ and $\cU_{\gamma_1} =
\cU_{\gamma_2}$. \hspace{6mm} & {\rm(C')} $\dim \cS_{\gamma_1} <
\dim \cS_{\gamma_2}$.
\\[0.5ex]
{\rm(D)} $\cS_\gamma = \cS_{\gamma_2}$ and $\cU_\gamma =
\cU_{\gamma_2}$ & {\rm(D')} $\dim \cU_{\gamma_1} > \dim
\cU_{\gamma_2}$.
\\
\hspace{6.5mm}for $\gamma \in [\gamma_1,\gamma_2]$. &
\end{tabular}\end{lemma}

Now we are in a position to state and prove our main theorem on the nonuniform dichotomy spectrum.

\begin{theorem}\label{main21}
The
\nds $\Sigma_{NED}(A)$ of {\rm\x{1.1}} is a disjoint union of $n$
closed intervals {\rm(}called {\it spectral intervals}{\rm)} where $0 \leq n
\leq N$, i.e., either $\Sigma_{NED}(A) = \emptyset$, or $\Sigma_{NED}(A)
=\R$, or $\Sigma_{NED}(A)$ is in one of the four cases
 \begin{equation*}\label{spectrumNED}
\Sigma_{NED}(A) = \left\{\begin{matrix}
[a_1,b_1]\\
\text{or}\\
(-\oo,b_1]\end{matrix}\right\}
\cup [a_2,b_2] \cup \cdots \cup [a_{n-1},b_{n-1}] \cup
\left\{\begin{matrix}
[a_n,b_n]\\
\text{or}\\
{[}a_n,\infty)\end{matrix}\right\},\end{equation*}
where $0<a_1 \leq b_1 < a_2 \leq b_2 < \cdots < a_n \leq b_n$.
Furthermore,
choose a
\be
\lb{bc1}
\gamma_0 \in \rho_{NED}(A)~ with~
(-\oo,\gamma_0)\subset \rho_{NED}(A);
\ee
otherwise define $\cU_{\gamma_0}:=\R \times \R^{N}$, $\cS_{\gamma_0}:=\R
\times \{0\}$,
and choose a
\be\lb{bc2}
\gamma_n \in \rho_{NED}(A)~ {
with}~ (\gamma_n,+\oo)\subset \rho_{NED}(A);
\ee
otherwise define $\cU_{\gamma_n}:=\R \times \{0\}$,
$\cS_{\gamma_0}:=\R \times \R^{N}$.
Then the sets
\[
\cW_0:=\cS_{\gamma_0}\quad{\rm and}\quad
\cW_{n+1}:=\cS_{\gamma_n}
\]
are both linear integral manifolds of {\rm \x{1.1}}. For $n\geq 2$,
choose $\gamma_i \in \rho_{NED}(A)$ with \be\lb{bc3}b_i < \gamma_i <
a_{i+1}\quad {for} \quad i=1,\ldots,n-1.
\ee
Then for every
$i=1,\ldots,n-1$ the intersection
\[
\cW_i:=\cU_{\gamma_{i-1}}\cap
\cS_{\gamma_{i}}
\]
is a linear integral manifold of {\rm \x{1.1}}
with $\dim \cW_i \geq 1$.
Moreover,
those linear integral manifolds $\cW_i, i=0,\ldots,n+1$,
called spectral manifolds,
are independent of the choice of $\gamma_0,\ldots,\gamma_n$ in
{\rm\x{bc1}}, {\rm\x{bc2}} and {\rm\x{bc3}}
and satisfy
\[
\cW_0\oplus\cdots\oplus\cW_{n+1}=\R \times \R^{N}
\]
in the sense of Whitney sum, i.e.,
$\cW_0+\cdots+\cW_{n+1}=\R\times\R^{N}$ but
$\cW_i \cap \cW_j=\R\times \{0\}$ for $i\neq j$.
\end{theorem}

\noindent\prf{Recall that the resolvent set $\rho_{NED}(A)$
is open and therefore $\Sigma_{NED}(A)$ is
the disjoint union of closed intervals. Next we will show that
$\Sigma_{NED}(A)$ consists of at most $N$ intervals. Indeed, if
$\Sigma_{NED}(A)$ contains $N+1$ components, then one can choose a
collection of points $\zeta_{1},\ldots,\zeta_{N}$ in
$\rho_{NED}(A)$ such that $\zeta_{1}< \cdots < \zeta_{N}$ and each
of the intervals $(-\oo,\zeta_{1}),(\zeta_{1},\zeta_{2}),\ldots,$
$(\zeta_{N-1},\zeta_{N}),(\zeta_{N},\infty)$ has nonempty
intersection with the spectrum $\Sigma_{NED}(A)$. Now Alternative II of Lemma \ref{lemma2.3}] implies \[0 \leq \dim \cS_{\zeta_1} < \cdots
< \dim \cS_{\zeta_N} \leq N\] and therefore either $\dim
\cS_{\zeta_1}=0$ or $\dim \cS_{\zeta_N}=N$ or both. Without loss of
generality, $\dim \cS_{\zeta_N}=N$, i.e., $\cS_{\zeta_N}=\R \times
\R^{N}$. Assume that \[\dot{x}=[A(t)-\zeta_N I]x\] admits a strong
\ned~ with invariant projector $P \equiv {\rm Id}$, then
\[\dot{x}=[A(t)-\zeta I]x\] also admits a \ned~ with the same
projector for every $\zeta>\zeta_N$. Now we have the conclusion
$(\zeta_{N},\infty)\subset \rho_{NED}(A)$, which is a contradiction.
This proves the alternatives for $\Sigma_{NED}(A)$.

Due to Lemma \ref{lemma2.3}, the sets $\cW_0,\ldots,\cW_{n+1}$ are
linear integral manifolds. To prove that $\dim \cW_1\geq
1,\ldots, \dim \cW_n\geq 1$ for $n \geq 1$, we assume that $\dim
\cW_1=0$, i.e., $\cU_{\gamma_0} \cap \cS_{\gamma_1}=\R \times
\{0\}$. If $(0,b_1]$ is a spectral interval this implies that
$\cS_{\gamma_1}=\R \times \{0\}$. The projector of the \ned~ of
\[\dot{x}=[A(t)-\gamma_1 I]x\] is $0$ and then we get
the contradiction $(-\oo,\gamma_1)\subset \rho_{NED}(A)$. If $[a_1,b_1]$
is a spectral interval then $[\gamma_0,\gamma_1]\cap \Sigma_{NED}(A)
\neq \emptyset$ and Alternative II of Lemma \ref{lemma2.3} yields a
contradiction. Therefore $\dim \cW_1\geq 1$ and similarly $\dim
\cW_n\geq 1$. Furthermore for $n\geq 3$ and $i=2,\ldots,n-1$ one has
$[\gamma_{i-1},\gamma_{i}]\cap \Sigma_{NED}(A) \neq \emptyset$ and
again Alternative II of Lemma \ref{lemma2.3} yields $\dim \cW_i\geq1$.

For $i<j$ we have $\cW_{i}\subset \cS_{\gamma_i}$ and
$\cW_{i}\subset \cU_{\gamma_{j-1}} \subset \cU_{\gamma_i}.$ Using Lemma \ref{lem2.1}, we have
 $\cW_{i} \cap \cW_{j}\subset
\cS_{\gamma_i}\cap \cU_{\gamma_i}=\R \times \{0\}$ and therefore $\cW_{i} \cap
\cW_{j} =\R \times \{0\}$ for $i \neq j$.

To show that $\cW_0\oplus\cdots\oplus\cW_{n+1}=\R \times \R^{N}$,
recall the monotonicity relations $\cS_{\gamma_0}\subset \cdots
\subset \cS_{\gamma_n}$,
$\cU_{\gamma_0}\supset\cdots\supset\cU_{\gamma_n}$, and the identity
$\cS_\gamma \oplus \cU_\gamma = \R \times \R^N$ for $\gamma \in
\R$. Therefore $\R \times \R^{N}=\cW_0 \times \cU_{\gamma_0}$.
Now we have\beaa
\R \times \R^{N} \EQ \cW_0 + \cU_{\gamma_0} \cap [\cS_{\gamma_1}+\cU_{\gamma_1}]\\
\EQ \cW_0 + [\cU_{\gamma_0} \cap \cS_{\gamma_1}]+\cU_{\gamma_1}\\
\EQ \cW_0 + \cW_1+\cU_{\gamma_1}. \eeaa
Doing the same for $\cU_{\gamma_1}$, we get
\beaa\R \times \R^{N} \EQ \cW_0 + \cW_1 + \cU_{\gamma_1} \cap [\cS_{\gamma_2}+\cU_{\gamma_2}]\\
\EQ \cW_0 + \cW_1 + [\cU_{\gamma_1} \cap \cS_{\gamma_2}]+\cU_{\gamma_2}\\
\EQ \cW_0 + \cW_1 + \cW_2+\cU_{\gamma_2}, \eeaa and mathematical
induction yields  $\R\times \R^{N}=\cW_0+\cdots+\cW_{n+1}$. To
finish the proof, let $\tilde{\gamma}_0,\ldots,\tilde{\gamma}_n \in
\rho_{NED}(A)$ be given with the properties  {\rm\x{bc1}},
{\rm\x{bc2}} and {\rm\x{bc3}}. Then Alternative I of Lemma
\ref{lemma2.3}] implies\[\cS_{\gamma_i}=\cS_{\tilde{\gamma}_i}
\quad {\rm and} \quad \cU_{\gamma_i}=\cU_{\tilde{\gamma}_i} \quad {\rm for} \quad i=0,\ldots,n\]
and therefore the linear integral manifolds $\cW_0,\ldots,\cW_{n+1}$
are independent of the choice of $\gamma_0,\ldots,\gamma_n$ in
{\rm\x{bc1}}, {\rm\x{bc2}} and {\rm\x{bc3}}.} $\Box$

\begin{definition} \label{Def3.2}
We say that {\rm(\ref{1.1})} is nonuniformly exponentially bounded if
there exist constants $K >0,\varepsilon\geq 0$ and $a \geq 0$ such that
\begin{equation}  \label{bounded-growth}
\| \Phi(t,s) \| \leq K e^{a |t-s|}e^{\varepsilon |s|},
\qquad \mbox{for } t,s \in  \mathbb{R}.
\end{equation}\end{definition}

\begin{lemma}\label{lemma3.2}
Assume that {\rm \x{1.1}} is nonuniformly exponentially bounded.
Then $\Sigma_{NED}(A)$ is a bounded closed set and $\Sigma_{NED}(A)\subset [-a,a]$.
\end{lemma}

\noindent{\bf Proof.}\, Assume that (\ref{bounded-growth}) holds. Let $\gamma > a$ and
$\alpha := \gamma - a > 0$, estimate
(\ref{bounded-growth}) implies
\[\| \Phi_\gamma(t,s) \| \leq K e^{-\alpha(t-s)}e^{\varepsilon |s|},
\quad \mbox{for } t \geq s\]
and therefore (2.8)$_\gamma$ admits a \ned~ with
invariant projector $P ={\rm Id}$. We have $\gamma \in \rho_{NED}(A)$ and
similarly for $\gamma < -a$, therefore $\Sigma_{NED}(A) \subset [-a,a]$. $\Box$

\begin{corollary}\label{BG}
Assume that {\rm\x{1.1}} is nonuniformly exponentially bounded. Then the nonuniform dichotomy spectrum $\Sigma_{NED}(A)$
of {\rm\x{1.1}} is the disjoint union of $n$ closed intervals
where $0 \leq n \leq N$, i.e.,
\[\Sigma_{NED}(A) =[a_1,b_1]\cup [a_2,b_2] \cup \cdots \cup [a_{n-1},b_{n-1}]
\cup[a_n,b_n],\]where $a_1 \leq b_1 < a_2 \leq b_2 < \cdots < a_n \leq b_n$.
\end{corollary}

Finally we present an example to illustrate that $\Sigma_{NED}(A)\neq\Sigma_{ED}(A)$ can occur.

\noindent {\bf Example 2.1}. Consider the scalar equation
$\dot{x}=A(t)x$ with $A(t)=\lambda_0+at\sin t$, where
$\lambda_0<a<0$ ($|a|\ll1$ is sufficiently small). Then
$\Sigma_{NED}(A)=[\lambda_0+a,\lambda_0-a]$ and $\Sigma_{ED}(A)=\R$.

In fact, the evolution operator of $\dot{x}=A(t)x$ is given by
\[\Phi(t,s)=
e^{\lambda_0 (t-s)- a\cos t(t-s)  - a s (\cos t-\cos s)+a(\sin t-\sin s)}.\]
For any $\gamma\in \R$, the evolution operator of the shifted system $\dot{x}=[A(t)-\gamma]x$ is given by
\begin{equation}\label{example1}
 \Phi_\gamma(t,s)=
e^{(-\gamma+\lambda_0) (t-s)- a \cos t(t-s)  - a s (\cos t-\cos s)+a (\sin t-\sin s)}.
\end{equation}
For any $\gamma\in (\lambda_0-a,+\infty)$,  it follows from (\ref{example1}) that
\[| \Phi_\gamma(t,s)|\leq e^{2|a|}e^{-(\gamma-\lambda_0+a)(t-s) }e^{2|a|\cdot |s|},\qquad t\geq s,\]
 which implies that the
 shifted system $\dot{x}=[A(t)-\gamma]x$ admits a \ned~ with $P=1$, by
 taking \[K=e^{2|a|},\quad \alpha=\gamma-\lambda_0+a>0,\quad \varepsilon=2|a|>0.\]
 Thus,\begin{equation}\label{P1}(\lambda_0-a,+\infty)\subset \rho_{NED}(A).
 \end{equation}

For any $\widetilde{\gamma}\in (-\infty, \lambda_0+a)$,  it follows from (\ref{example1}) that
\[| \Phi_{\widetilde{\gamma}}(t,s)|\leq e^{2|a|}e^{(-\widetilde{\gamma}+\lambda_0+a)(t-s) }e^{2|a|\cdot |s|},\,\,\,\,\mbox{for}\,\,\, t\leq s,\]
which implies that the shifted system $\dot{x}=[A(t)-\widetilde{\gamma}]x$ admits a \ned~ with $P=0$,
 by taking \[K=e^{2|a|},\quad \widetilde{\alpha}=-\widetilde{\gamma}+\lambda_0+a>0, \quad\varepsilon=2|a|>0.\]
Thus,
\begin{equation}\label{P0}(-\infty, \lambda_0+a)\subset \rho_{NED}(A).\end{equation}
It follows from (\ref{P1})-(\ref{P0}) that
\[(-\infty, \lambda_0+a) \cup   (\lambda_0-a,+\infty)\subset \rho_{NED}(A),\]
which implies that
\[\Sigma_{NED}(A) \subset[ \lambda_0+a,  \lambda_0-a].\] Now we show that
\[[ \lambda_0+a,  \lambda_0-a] \subset \Sigma_{NED}(A).\]
To show this, we first prove that $ \lambda_0-a  \in \Sigma_{NED}(A)$. On the contrary, assume
 that $\gamma_2=  \lambda_0-a$ such that
$\dot{x}=[A(t)-\gamma_2]x$ admits a \ned. We know that either the projector
$P=0$ or $P=1$. If $P=1$, then there exist constants $K,\,\alpha>0$ and $\varepsilon>0$ such that the following estimate holds
\beaa|\Phi_{\gamma_2}(t,s)|\EQ e^{[-\gamma_2+\lambda_0] (t-s)- a \cos t(t-s)  -
a s (\cos t-\cos s)+a (\sin t-\sin s)}\\
\LE K e^{-\alpha (t-s)} e^{\varepsilon |s|},\quad t\geq s.\eeaa
Substituting $\gamma_2=\lambda_0-a$, we have
\[e^{ a(1- \cos t)(t-s)  - a s (\cos t-\cos s)+a (\sin t-\sin s)}\leq K e^{-\alpha (t-s)} e^{\varepsilon |s|},
\,\,\, t\geq s,\]
which yields a contradiction for $s=0$ and $t\rightarrow +\infty$. If $P=0$, the dichotomy estimate is
\[e^{ a(1- \cos t)(t-s)  - a s (\cos t-\cos s)+a (\sin t-\sin s)}\leq K e^{\alpha (t-s)} e^{\varepsilon |s|},
\,\,\, t\leq s,\]
which also yields a contradiction for $s=-(2k-1)\pi$ and $t=-2k\pi$ and $k\rightarrow +\infty$.
Therefore, $ \lambda_0-a  \in \Sigma_{NED}(A)$. Analogously, we can
prove that $ \lambda_0+a  \in \Sigma_{NED}(A)$. By Theorem \ref{main21},
we know that $\Sigma_{NED}(A)$ is an interval.
Thus, for any $\gamma\in [ \lambda_0+a,  \lambda_0-a]$, it follows from
the connectedness that $ \gamma\in\Sigma_{NED}(A)$. Consequently,  \[
[ \lambda_0+a,  \lambda_0-a] \subset \Sigma_{NED}(A).\]
Therefore, $\Sigma_{NED}(A)=[ \lambda_0+a,  \lambda_0-a]$.

On the other hand, we can show that, for any $\gamma\in (\lambda_0-a,+\infty)\cup(-\infty,\lambda_0+a)$,
the shifted system $\dot{x}=[A(t)-\gamma]x$
admits no exponential dichotomy. From the above proof, $\Sigma_{NED}(A)=[ \lambda_0+a,  \lambda_0-a]$, which
implies that the shifted system $\dot{x}=[A(t)-\gamma]x$ admits no \ned.
Consequently, for $\gamma\in [ \lambda_0+a,  \lambda_0-a]$,  the shifted
system $\dot{x}=[A(t)-\gamma]x$ admits no exponential dichotomy.
Therefore, $\Sigma_{ED}(A)=\R$.


\section{Reducibility}
\setcounter{equation}{0}

In this section we employ Theorem \ref{main21} to prove a reducibility result. We refer to \cite{cop67, pal, s-jlms} and
the references therein for some reducibility results in the
setting of classic exponential dichotomies.
First we recall the definition of {\em kinematic similarity} and
several results in \cite{s-jlms}.

\begin{definition} \label{Def1.2} Given $A,B\in \mathfrak{L}_{loc}^1$.
Equation {\rm(\ref{1.1})} is said to be kinematically similar to another system
\be y'=B(t)y,\label{1.3}\ee
if there exists an absolutely continuous function  $S:\mathbb{R}\rightarrow GL_N(\mathbb{R})$
with $\sup_{t\in \mathbb{R}}\|S(t)\|<\infty$ and $\sup_{t\in \mathbb{R}}\|S^{-1}(t)\|<\infty$ which satisfies the
differential equation
\be\lb{SAB}S'(t)=A(t)S-SB(t).\ee
The transformation $x=S(t)y$ which transforms {\rm(\ref{1.1})} into {\rm(\ref{1.3})} is called the Lyapunov transformation.
\end{definition}

\begin{lemma}\label{lemma1.2} \cite[Lemma A.5]{s-jlms} Let $P\in \mathbb{R}^{N\times N}$
be a symmetric projection and $X: \mathbb{R}\rightarrow GL_N(\mathbb{R})$ be an absolutely
continuous matrix. Then
\bb{itemize} \item[{\rm(A)}] The mapping
\[\widetilde{R}:\R\rightarrow\R^{N \times N},\quad t \mapsto PX(t)^T X(t) P +QX(t)^T X(t)Q\]
is absolutely continuous and $\widetilde{R}(t)$ is a positive definite, symmetric
matrix for every $t \in \mathbb{R}$. Moreover there is a unique absolutely
continuous function
$R : \mathbb{R}\rightarrow \mathbb{R}^{N \times N}$
of positive definite symmetric matrices $R(t)$, $t \in \mathbb{R}$, with
\[R(t)^2 = \widetilde{R}(t),\quad PR(t) = R(t) P.\]
\item[{\rm(B)}]The mapping
\[S:\mathbb{R}\rightarrow \mathbb{R}^{N\times N},
\quad t \mapsto X(t) R(t)^{-1}\]
is absolutely continuous and $S(t)$ is invertible, satisfying
\[S(t)PS^{-1}(t)=X(t)PX^{-1}(t),\]
\[S(t)QS^{-1}(t)=X(t)QX^{-1}(t),\]
\[\|S(t)\|{\leq}\sqrt{2},\]
\[\|S^{-1}(t)\|{\leq}[\|X(t)PX^{-1}(t)\|^{2}+\|X(t)QX^{-1}(t)\|^{2}]^{\frac{1}{2}},\quad t\in\R.\]
\end{itemize}
\end{lemma}

In the setting of classical exponential dichotomies, $S^{-1}(t)$ is bounded, which
follows from the properties
$\|X(t)PX^{-1}(t)\|<\oo$ and $\|X(t)QX^{-1}(t)\|<\oo$.
However, in the setting of nonuniform exponential dichotomies, $S^{-1}(t)$ can be unbounded, because
\[\|X(t)PX^{-1}(t)\|\leq Ke^{\varepsilon t},\quad t\geq 0.\]
Such a fact will make difficulties to the analysis. To overcome it, we introduce
the new notions of {\em nonuniform  Lyapunov transformation} and
{\em nonuniform kinematical similarity}.

\begin{definition} \label{Def2.1}
Suppose that  $S:\mathbb{R}\rightarrow GL_N(\mathbb{R})$ is an absolutely
continuous matrix. $S(t)$ is said to be a nonuniform Lyapunov matrix if
there exists a constant $M=M_{\varepsilon}>0$ such that
\[\|S(t)\|\leq Me^{\varepsilon |t|}\quad and\quad  \|S^{-1}(t)\| \leq Me^{\varepsilon |t|},
\quad\mbox{for all}~t\in\R.\]
\end{definition}

\begin{definition} \label{Def2.2} Equation {\rm(\ref{1.1})} is said to be nonuniformly
kinematically similar to equation {\rm(\ref{1.3})} if there exists a nonuniform Lyapunov matrix $S(t)$
satisfying the differential equation {\rm\x{SAB}}.
For short, we write {\rm(\ref{1.1})} $\sim$ {\rm(\ref{1.3})} or $A(t) \sim B(t)$.
\end{definition}
For the sake of comparison, we denote kinematical similarity by {\rm(\ref{1.1})} $\approx$ {\rm(\ref{1.3})} or $A(t) \approx B(t)$.

\begin{definition} \label{Def2.3}
We say that system {\rm(\ref{1.1})} is reducible, if it is nonuniformly kinematically
similar to system {\rm(\ref{1.3})} whose coefficient matrix $B(t)$ has the block form
\begin{equation}\left(\begin{array}{lll}
B_1(t) & 0\\0 &  B_2(t)
\end{array}\right),\label{2.1}\end{equation}
where $B_1(t)$ and $B_2(t)$ are both matrices of smaller size than $B(t)$.
\end{definition}

In \cite{cop78}, Coppel proved that if system (\ref{1.1}) admits an
exponential dichotomy, then there exists a Lyapunov transformation such
that $A(t) \approx B(t)$ and $B(t)$ has the block form (\ref{2.1}),
i.e., system (\ref{1.1}) is {\em reducible}.
The following theorem shows that if system (\ref{1.1}) admits a nonuniform
exponential dichotomy, then there exists a nonuniform Lyapunov transformation
such that  $A(t)\sim B(t)$ and $B(t)$ has the block form (\ref{2.1}),
i.e., system (\ref{1.1}) is reducible.

\begin{lemma}\label{lemma3.4}  Suppose that system {\rm(\ref{1.1})} admits a strong
\ned~ with the form of estimates {\rm\x{np}-\x{nq}}
and $rank(\tilde{P})=k,(0\leq k\leq N)$. If there exists a nonuniform Lyapunov
transformation $S(t)$ such that $A(t)\sim B(t)$, then system {\rm(\ref{1.3})} also
admits a \ned~and the projector has the same rank.
\end{lemma}

\noindent{\bf Proof.} Suppose that $S(t)$ is the nonuniform  Lyapunov matrix with $\|S(t)\|\leq Me^{\varepsilon |t|},\|S^{-1}(t)\|  \leq Me^{\varepsilon |t|}$ and such that $A(t)\sim B(t)$.
Let $Y(t)=S(t)X(t)$. Then it is easy to see that $Y(t)$ is the fundamental matrix of
system (\ref{1.3}). To prove that system (\ref{1.3}) admits a \ned, we first consider the case $t\geq 0$. For $t\geq 0$,
\begin{equation}\begin{array}{lll}
\|Y(t)\tilde{P}Y^{-1}(s)\|&=&
\|S(t)X(t)\tilde{P}X^{-1}(s)S^{-1}(s)\|\\
&\leq&\|S(t)\|\cdot\|X(t)\tilde{P}X^{-1}(s)\|\cdot\|S^{-1}(s)\|\\
&\leq& K M^2 e^{\varepsilon |t|} e^{-\alpha(t-s)}e^{\varepsilon |s|}e^{\varepsilon |s|}\\
&\leq& K M^2 e^{\varepsilon (t-s)} e^{-\alpha(t-s)}e^{3\varepsilon |s|}\\
&=& K M^2 e^{-(\alpha-\varepsilon)(t-s)}e^{ 3\varepsilon    |s|},\,\,t\geq s.
\end{array}\label{2.5}\end{equation}
A similar argument shows that
\begin{equation}
\|Y(t)\tilde{Q}Y^{-1}(s)\|\leq  K M^2 e^{(\alpha+\varepsilon)(t-s)}e^{3\varepsilon  |s|},\,\,t\leq s.
\label{2.6}\end{equation}
It follows from (\ref{2.5})-(\ref{2.6}) that (\ref{1.3}) admits a \ned~
for $t\geq0$ due to $\varepsilon<\alpha$. Similarly, we see that system (\ref{1.3}) admits a \ned~ for $t\leq0$.
Thus (\ref{1.3}) admits a \ned~ and the rank of the projector is $k$. $\Box$

\begin{corollary}\label{cor3.2}
Assume that there exists a nonuniform  Lyapunov transformation $S(t)$
such that $A(t)\sim B(t)$. Then $\Sigma_{NED}(A) =  \Sigma_{NED}(B)$.
\end{corollary}

\begin{theorem}\label{Th2.1} Assume that equation  {\rm(\ref{1.1})} admits a nonuniform exponential dichotomy of the
 form {\rm\x{np}-\x{nq}} with invariant projector $P(t)\neq0,I$. Then {\rm(\ref{1.1})} is nonuniformly kinematically
similar to a decoupled system
\[\dot{x} =\begin{pmatrix}
B_1(t) & 0 \\0 & B_2(t)
\end{pmatrix}x\]
for some locally integrable matrix functions
\[B_1 :  \mathbb{R} \rightarrow  \mathbb{R}^{N_1 \times N_1}
\quad \text{and} \quad
B_2 :  \mathbb{R} \rightarrow  \mathbb{R}^{N_2 \times N_2}\]
where $N_1 := \dim \im P$ and $N_2 := \dim \ker P$.
\end{theorem}

\noindent{\bf Proof.} Since equation  (\ref{1.1}) admits a \ned{} of the form \x{np}-\x{nq}
with invariant projector $P(t)\neq0,I$, by Lemma \ref{lemma2.1}, we can choose a fundamental matrix $X(t)$
and the projector $P_0= \left(\begin{array}{lll}
I_{N_1} & 0\\
0 &  0\end{array}\right)$, $(0<k<N)$ such that the estimates \x{xnp}-\x{xnq} hold.
For the given nonsigular matrix $X(t)$, by Lemma \ref{lemma1.2}, there exists an
absolutely continuous and invertible matrix $S(t)$ satisfying
\[\|S(t)\|{\leq}\sqrt{2} ,\]
\[\|S^{-1}(t)\|{\leq}[\|X(t)\tilde{P}X^{-1}(t)\|^{2}+\|X(t)\tilde{Q}X^{-1}(t)\|^{2}]^{\frac{1}{2}} ,\]
which combined with the estimates \x{ntp}-\x{ntq} gives
\[\|S(t)\|{\leq}\sqrt{2}\leq M e^{\varepsilon |t|},\]
\[\|S^{-1}(t)\|{\leq}[\|X(t)\tilde{P}X^{-1}(t)\|^{2}+\|X(t)\tilde{Q}X^{-1}(t)\|^{2}]^{\frac{1}{2}}\leq \sqrt{2}Ke^{\varepsilon |t|}.\]
Thus we can take $M=M_\varepsilon\geq \max\{\sqrt{2},\sqrt{2}K\}$ such that
\[\|S(t)\|{\leq}  M e^{\varepsilon |t|},\quad \|S^{-1}(t)\|{\leq}  M e^{\varepsilon |t|},\]
which implies that $S(t)$ is a nonuniform  Lyapunov matrix. Setting
 \be\lb{br}B(t)=\dot{R}(t)R^{-1}(t), \ee
 where $R(t)=S(t)X(t)$ and define $B(t)=0$ for $t\in \mathbb{R}$ for which $\dot{S}(t)$ does not
exist. Obviously, $R(t)$ is the fundamental matrix of the linear equation
\[\dot{y}=B(t)y.\]
Now we show that $A(t)\sim B(t)$ and $B(t)$ has the block diagonal form
\[B(t)= \left(\begin{array}{lll}B_1(t) & 0\\
0 &  B_{2}(t)\end{array}\right),\quad t\in\R.\]
In fact,
\[\begin{array}{lll}
S'(t)&=&(X(t)R^{-1}(t))'\\
&=& X'(t)R^{-1}(t)+X(t)(R^{-1}(t))'\\
&=&A(t)X(t)R^{-1}(t)-X(t)R^{-1}(t)R'(t)R^{-1}(t),
\end{array}\]
which, combining with \x{br} gives
\[S'(t)=A(t)S(t)-S(t)B(t).\]
Therefore, $A(t)\sim B(t)$. Now we show that $B(t)$ has the block diagonal form of \x{2.1}. By Lemma \ref{lemma1.2}, $R(t)$
and $R(t)^{-1}$ commute with the matrix $\tilde{P}$ for every $t \in\R$.
The derivatives $\dot{R}(t)$ also commute with $\tilde{P}$, and then
\be\lb{PB}\tilde{P}B(t) = B(t)\tilde{P}\ee for almost all $t\in\mathbb{R}$.
Now we decompose
$B: \R \rightarrow \R^{N \times N}$ into four functions
\[\begin{array}{lll}
B_1: \R\rightarrow \R^{N_{1} \times N_{1}}, & B_2: \R
\rightarrow \R^{N_{2} \times N_{2}},\\
B_3: \R \rightarrow \R^{N_{1} \times N_{2}}, & B_4: \R
\rightarrow \R^{N_{2} \times N_{1}},
\end{array}\]
with\[B(t)= \left(\begin{array}{lll} B_1(t) & B_3(t)
\\B_4(t) &  B_2(t)\end{array}\right), \quad t\in\R.\]
Identity \x{PB} implies that
\[\left(\begin{array}{lll} B_1(t) & B_3(t)\\
0 &  0\end{array}\right)=\left(\begin{array}{lll} B_1(t) & 0
\\B_4(t) &  0\end{array} \right), \quad t\in\R.\]
Therefore $B_3(t)\equiv 0$ and $B_4(t)\equiv 0$. Thus $B$ has
the block diagonal form
\[B_{k}= \left(\begin{array}{lll} B_1(t) & 0\\
0 &  B_2(t)\end{array} \right), \quad t\in\R\]
and the proof is finished. $\Box$

Now we are in a position to prove the reducibility result.
\begin{theorem} Assume that {\rm(\ref{1.1})} admits a nonuniform
exponential dichotomy. Due to Theorem {\rm
\ref{main21}}, the dichotomy spectrum is either empty or the
disjoint union of $n$ closed spectral intervals $\cI_1, \dots,
\cI_n$ with $1 \leq n \leq N$, i.e.,
\[\Sigma_{NED}(A) = \emptyset \quad(n = 0)
\qquad\text{ or }\qquad
\Sigma_{NED}(A) = \cI_1 \cup \cdots \cup \cI_n\;.\]
Then there exists a weakly kinematic
similarity action $S : \R \to \R^{N \times N}$ between
{\rm(\ref{1.1})} and a block diagonal system
\[ \label{nsystem}\dot{x}=\begin{pmatrix}
B_{0}(t) & & \\& \ddots & \\
 & & B_{n+1}(t)\end{pmatrix}x\]with locally integrable functions $B_i: \R\rightarrow
\R^{N_i \times N_i}$, $N_i=\dim \cW_i$, and
\[\Sigma_{NED}(B_0) = \emptyset,\quad\Sigma_{NED}(B_1) = \cI_1, \cdots,
\Sigma_{NED}(B_n) =\cI_n,\quad\Sigma_{NED}(B_{n+1}) = \emptyset.\]
\end{theorem}

\noindent\prf{If for any $\gamma\in\R$,  system (2.8)$_\gamma$ admits a
\ned, then $\Sigma_{NED}(A) = \emptyset$. Conversely, for any
$\gamma\in \R$, the weighted system (2.8)$_\gamma$ does not admit a \ned,
then $\Sigma_{NED}(A) =\R$. Now, we prove the theorem for the
nontrivial case ($\Sigma_{NED}(A) \neq \emptyset$ and
$\Sigma_{NED}(A) \neq\R$). First, recall that the resolvent set
$\rho_{NED}(A)$ is open and therefore the dichotomy spectrum
$\Sigma_{NED}(A)$ is the disjoint union of closed intervals. Using
Theorem \ref{main21}, we can assume\[\cI_{1} = \left\{\begin{matrix}
    [a_1,b_1]
    \\
    \text{or}
    \\
    (-\oo,b_1]\end{matrix}\right\}
, \cI_{2}=[a_2,b_2], \cdots , \cI_{n-1}=[a_{n-1},b_{n-1}] ,
  \cI_{n} =\left\{\begin{matrix}
    [a_n,b_n]
    \\
    \text{or}
    \\
    {[}a_n,\infty)\end{matrix}\right\}\]
with $a_1 \leq b_1 < a_2 \leq b_2 < \cdots < a_n \leq b_n$ due to
$\varepsilon<\alpha/2$.

If $\cI_{1}=[a_1,b_1]$ is a spectral interval, then we have$(-\oo,
\gamma_0)\subset \rho_{NED}(A)$ and $\cW_0=\cS_{\gamma_0}$ for some
$\gamma_0 < a_1$ due to Theorem \ref{main21}, which implies that
\[\dot{x} = [A(t) - \gamma_{0} I]x\] admits a \ned ~
with an invariant projector $\tilde{P}_0$. By Corollary \ref{cor3.2}
and Theorem \ref{Th2.1}, there exists a  nonuniform Lyapunov
transformation $x=S_{0}(t)x_{0}$ with $\|S_{0}(t)\|\leq
M_{0}e^{\epsilon_{0}|t|}$ and $\|S_{0}(t)^{-1}\| \leq
M_{0}e^{\epsilon_{0}|t|}$  such that $A(t) \sim A_{0}(t)$ and
$A_{0}(t)$ has two blocks of the form
$A_{0}(t)=\left(\begin{array}{lll}
 B_{0}(t) & 0
 \\
 0 &   B_{0,*}(t)
 \end{array}\right)
$ with $\dim B_{0}(t)=\dim \im \tilde{P}_0=\dim \cS_{\gamma_0}=\dim
\cW_{0}:=N_{0}$ due to Theorem~\ref{Th2.1}, Lemma~\ref{lem2.1} and
Theorem \ref{main21}. If $\cI_{1}=(-\oo,b_1]$ is a spectral
interval, a block $B_{0}(t)$ is omitted.

Now we consider the following system
\[\dot{x}_{0} = A_{0}(t)x_{0}=\left(\begin{array}{lll}
 B_{0}(t) & 0
 \\
 0 &   B_{0,*}(t)
 \end{array}\right)x_{0}.\]
By using Lemma \ref{lemma2.3}, we take $\gamma_1\in (b_{1},a_{2})$.
In view of $(b_{1},a_{2})\subset \rho_{NED}(B_{0,*}(t))$,
$\gamma_1\in \rho_{NED}(B_{0,*}(t))$, which implies that
\[\dot{x}_{0} = \left[\left(\begin{array}{lll}
 B_{0}(t) & 0
 \\
 0 &   B_{0,*}(t)
 \end{array}\right) - \gamma_{0} I\right]x_{0}\]
admits a \ned~ with an invariant projector $\tilde{P}_1$. From the
claim above, we know that $\tilde{P}_1\neq 0,\,I$. Similarly by
Corollary \ref{cor3.2} and Theorem \ref{Th2.1}, there exists a
nonuniform Lyapunov transformation
\[x_{0}=S_{1}(t)x_{1}=\left(\begin{array}{lll}
 I_{N_{0}} & 0
 \\
 0 & \tilde{S}_{1}(t)
 \end{array}
 \right)x_{1}\]
with $\|\tilde{S}_{1}(t)\|\leq M_{1}e^{\epsilon_{1}|t|}$ and
$\|(\tilde{S}_{1}(t))^{-1}\| \leq M_{1}e^{\epsilon_{1}|t|}$ such
that $ B_{0,*}(t)\sim \tilde{B}_{0,*}(t)$ and $\tilde{B}_{0,*}(t)$
has two blocks of the form
$\tilde{B}_{0,*}(t)=\left(\begin{array}{lll}
 B_{1}(t) & 0
 \\
 0 &   B_{1,*}(t)
 \end{array}
 \right)
$ with $B_{1}(t)=\dim \im \tilde{P}_1=\dim
\cS_{\gamma_1}\geq\dim(\cU_{\gamma_0}\cap \cS_{\gamma_1})=\dim
\cW_{1}:=N_{1}$ due to Theorem~\ref{Th2.1}, Lemma~\ref{lem2.1} and
Theorem \ref{main21}. In addition, using Corollary \ref{cor3.2} and
Theorem \ref{Th2.1}, we have
\[\Sigma_{NED}(B_{1}(t)) =\left\{\begin{matrix}
    [a_1,b_1]
    \\
    \text{or}
    \\
    (-\oo,b_1]
  \end{matrix}\right\},~~
\Sigma_{NED}(B_{1,*}(t)) =[a_2,b_2] \cup \cdots \cup
[a_{n-1},b_{n-1}] \cup
  \left\{\begin{matrix}
    [a_n,b_n]
    \\
    \text{or}
    \\
    {[}a_n,\infty)
  \end{matrix}\right\}.\]Now we can construct a nonuniform
Lyapunov transformation $x=\tilde{S}(t)x_{1}$ with $
 \tilde{S}(t)=S_{0}(t)S_{1}(t)=S_{0}(t)\left(\begin{array}{lll}
 I_{N_{0}} & 0
 \\
 0 & \tilde{S}_{1}(t)
 \end{array}
 \right)
$, where $\|\tilde{S}(t)\|\leq M_0
M_1e^{(\epsilon_{0}+\epsilon_{1})|t|}$ and $\| \tilde{S}(t)^{-1}\|
\leq M_0 M_1e^{(\epsilon_{0}+\epsilon_{1})|t|}$. Then $A(t)\sim
A_{1}(t)$ and $A_{1}(t)$ has three blocks of the form \[ A_{1}(t)=
\begin{pmatrix}
    B_{0}(t) & & \\
    & B_{1}(t) & \\
    & & B_{1,*}(t)
\end{pmatrix}.\]

Applying similar procedures to $\gamma_{2}\in(b_{2},a_{3})$,
$\gamma_{3}\in(b_{3},a_{4}),\ldots$, we can construct a weakly
non-degenerate transformation $x=S(t)x_n$ with \[
 S(t)=S_{0}(t)\left(\begin{array}{lll}
 I_{N_{0}} & 0
 \\
 0 & \tilde{S}_{1}(t)
 \end{array}
 \right)\left(\begin{array}{lll}
 I_{N_{0}+N_{1}} & 0
 \\
 0 & \tilde{S}_{2}(t)
 \end{array}
 \right)\cdots \left(\begin{array}{lll}
 I_{N_{0}+\ldots+N_{n-1}} & 0
 \\
 0 & \tilde{S}_{n}(t)
 \end{array}
 \right)\]
such that $\| S(t)\|\leq M e^{\epsilon|t|}$ and $\| S(t)^{-1}\| \leq
M e^{\epsilon|t|}$ with $M=M_0\times \cdots \times M_n$ and
$\epsilon=\epsilon_0+ \cdots+\epsilon_n$ . Now we can prove
\[A(t)\sim A_{n}(t):=B(t) =\begin{pmatrix}
B_{0}(t) & & \\
& \ddots & \\
& & B_{n+1}(t)\end{pmatrix}\] with locally integrable functions
$B_{i}: \R \rightarrow \R^{N_i \times N_i}$ and\[\Sigma_{NED}(B_0) =
\emptyset \;,\Sigma_{NED}(B_1) = \cI_1 \;, \dots \;,
\Sigma_{NED}(B_n) =\cI_n,\Sigma_{NED}(B_{n+1}) = \emptyset.\]

Finally, we show that $N_i=\dim \cW_i$. From the claim above, we
note that $\dim B_{0}(t)=\dim \cW_0, ~\dim B_{1}(t)\geq \dim
\cW_1,\ldots, \dim B_{n}(t)\geq \dim \cW_n, ~\dim B_{n+1}(t)=\dim
\cW_{n+1}$ and with Theorem \ref{main21} this gives $\dim \cW_0 +
\cdots +\dim \cW_{n+1}=N$, so $\dim B^{i}_k=\dim \cW_i$ for
$i=0,\ldots,n+1$. Now the proof is finished.}  $\Box$

\end{document}